\documentclass[11pt]{article}
\usepackage[utf8]{inputenc}
\usepackage{authblk}
\usepackage{mathtools}
\usepackage[english]{babel}
\usepackage{amsmath}
\usepackage{amsthm}
\usepackage{amsfonts}
\usepackage{amssymb}
\usepackage{bbm}
\usepackage{subcaption}
\usepackage{verbatim}
\usepackage[usenames,dvipsnames]{color}
\usepackage{hyperref}
\usepackage[left=2cm,right=2cm,top=2cm,bottom=2cm]{geometry}
\usepackage{dsfont}
\usepackage{enumitem}
\usepackage{tikz}
\usepackage[numeric,initials,nobysame,msc-links,abbrev]{amsrefs}
\renewcommand{\eprint}[1]{\href{https://arxiv.org/abs/#1}{arXiv:#1}}
\newcommand{\pageafter}[1]{#1~pp.}
\BibSpec{article}{%
+{} {\PrintAuthors} {author}
+{,} { \textit} {title}
+{.} { } {part}
+{:} { \textit} {subtitle}
+{,} { \PrintContributions} {contribution}
+{.} { \PrintPartials} {partial}
+{,} { } {journal}
+{} { \textbf} {volume}
+{} { \PrintDatePV} {date}
+{,} { \issuetext} {number}
+{,} { \pageafter} {pages}
+{,} { } {status}
+{,} { \PrintDOI} {doi}
+{,} { available at \eprint} {eprint}
+{} { \parenthesize} {language}
+{} { \PrintTranslation} {translation}
+{;} { \PrintReprint} {reprint}
+{.} { } {note}
+{.} {} {transition}
+{} {\SentenceSpace \PrintReviews} {review}
}
\BibSpec{collection.article}{%
+{} {\PrintAuthors} {author}
+{,} { \textit} {title}
+{.} { } {part}
+{:} { \textit} {subtitle}
+{,} { \PrintContributions} {contribution}
+{,} { \PrintConference} {conference}
+{} {\PrintBook} {book}
+{,} { } {booktitle}
+{,} { \PrintDateB} {date}
+{,} { \pageafter} {pages}
+{,} { } {status}
+{,} { \PrintDOI} {doi}
+{,} { available at \eprint} {eprint}
+{} { \parenthesize} {language}
+{} { \PrintTranslation} {translation}
+{;} { \PrintReprint} {reprint}
+{.} { } {note}
+{.} {} {transition}
+{} {\SentenceSpace \PrintReviews} {review}
}

\usetikzlibrary{arrows}
\usepackage[normalem]{ulem}

\newtheorem{theorem}{Theorem}

\newtheorem{proposition}[theorem]{Proposition}

\theoremstyle{definition}

\newtheorem{remark}[theorem]{Remark}

\newcommand{\bbB}{\mathbb{B}}

\newcommand{\bbN}{\mathbb{N}}

\newcommand{\bbP}{\mathbb{P}}

\newcommand{\bbR}{\mathbb{R}}

\newcommand{\bbW}{\mathbb{W}}

\newcommand{\bbZ}{\mathbb{Z}}

\newcommand{\cC}{\mathcal{C}}

\newcommand{\cE}{\mathcal{E}}

\newcommand{\cI}{\mathcal{I}}

\newcommand{\cP}{\mathcal{P}}

\newcommand{\cS}{\mathcal{S}}
\newcommand{\cT}{\mathcal{T}}

\newcommand{\1}{\mathbbm{1}}

\newcommand{\md}{\mathrm{d}}
\newcommand{\mS}{\mathrm{S}}
\newcommand{\mI}{\mathrm{I}}
\newcommand{\mR}{\mathrm{R}}

\DeclareMathOperator{\supp}{supp}

\newcommand{\eps}{\varepsilon}

\title{\scshape Brownian snails with removal die out in one dimension}

\author[1]{Ivailo Hartarsky}
\author[2,3]{Lyuben Lichev}
    \affil[1]{TU Wien, Faculty of Mathematics and Geoinformation, Institute of Statistics and Mathematical Methods in Economics, Research Unit of Mathematical Stochastics, Wiedner Hauptstra\ss e 8-10, A-1040 Vienna, Austria, \texttt{ivailo.hartarsky@tuwien.ac.at}}
\affil[2]{Institut Camille Jordan, Univ. Jean Monnet, Saint-Etienne, France, \texttt{lyuben.lichev@univ-st-etienne.fr}}
\affil[3]{Institute of Mathematics and Informatics, Bulgarian Academy of Sciences, Acad.~G.~Bonchev~str. 8, 1113 Sofia, Bulgaria}

\begin{document}

\maketitle

\begin{abstract}
Brownian snails with removal is a spatial epidemic model defined as follows. Initially, a homogeneous Poisson process of susceptible particles on $\bbR^d$ with intensity $\lambda>0$ is deposited and a single infected one is added at the origin. Each particle performs an independent standard Brownian motion. Each susceptible particle is infected immediately when it is within distance 1 from an infected particle. Each infected particle is removed at rate $\alpha>0$, and removed particles remain such forever. Answering a question of Grimmett and Li, we prove that in one dimension, for all values of $\lambda$ and $\alpha$, the infection almost surely dies out.
\end{abstract}


\section{Introduction}
\label{sec:intro}

Susceptible/infected/removed (SIR) models are among the most classical in epidemiology. While they are most commonly studied in a mean-field setting, considerable progress in the analysis of geometric SIR models for stochastic dynamics witnessing the individuals in the population was made in the past decades. Given the circumstances (parameter values), the question of utmost importance is whether the disease almost surely dies out or it survives indefinitely with positive probability.

The present paper deals with the Brownian snails with removal model recently studied by Grimmett and Li \cite{Grimmett22}. Before introducing the model formally, let us briefly discuss related simpler models and results. The \emph{frog model} \cites{Alves02,Alves02a,Popov03} is a lattice SIR system in which infected particles perform random walks (in discrete or continuous time) on $\bbZ^d$ until they are removed at rate $\alpha$ (in discrete or continuous time). 
Susceptible particles remain immobile until an infected particle visits the site they occupy, at which point they become infected (and start their own random walk until they are removed). We refer to the latter feature (susceptible particles not moving) as \emph{delay}. 
While it appears harmless from a modelisation viewpoint, considering a delayed model makes it mathematically easier to handle. Similarly, there is a canonically associated SI model (without removal), obtained by taking $\alpha=0$. Such models enjoy crucial monotonicity properties making them easier to study. 
For instance, for the frog model without removal, it is known that the set of infected particles converges to a (linearly expanding) limit shape \cites{Alves02,Ramirez02}, much like classical static models such as first passage percolation. 

Similarly to the frog model without removal, a delayed model without removal in continuous space (where particles perform Brownian motions rather than random walks and infections are transmitted within distance 1) was studied much more recently by Beckman, Dinan, Durrett, Huo and Junge \cite{Beckman18}, where it was proved that the set of infected particles also exhibits a limit shape. Grimmett and Li \cite{Grimmett22} investigated the corresponding model of delayed Brownian snails with removal, showing that it exhibits a non-trivial phase transition only for dimensions 2 and higher.

We next turn to models without delay. Among the first results in this direction are the ones of Kesten and Sidoravicius \cites{Kesten05,Kesten06,Kesten08} who proved a shape theorem in the absence of removal in a discrete space model. Gracar and Stauffer \cite{Gracar18} treated this model also in finite volume. The most recent progress is due to Dauvergne and Sly \cite{Dauvergne23} who showed that infection progresses linearly (though without establishing a limit shape result) also if one allows for different speeds for infected and susceptible particles, which creates major difficulties.
Subsequently, they also treated the corresponding model with removal for small enough removal rate $\alpha$ in dimensions $d\ge 2$ \cite{Dauvergne22}. In the setting of Brownian snails (with removal and no delay) Grimmett and Li \cite{Grimmett22} only treated the case $\alpha$ large enough, again in dimension $d\ge 2$.

However, determining whether Brownian snails have a non-trivial phase transition in one dimension remained beyond the scope of \cite{Grimmett22}. Indeed, in \cite{Grimmett22}*{Question D} they asked for a solution of this problem, which was subsequently reiterated in \cite{Grimmett23}*{Question 6.1(i)}. Our main result resolves this question by showing that in one dimension there is no survival regardless of the value of the removal rate $\alpha>0$ and the intensity $\lambda$ of the initial configuration.

Analogues of our result have been proved for the frog model with removal \cite{Alves02a}*{Theorem 1.1} and the delayed Brownian snails model with removal \cite{Grimmett22}*{Theorem 1.2} but not for non-delayed models. In \cite{Dauvergne22}*{Remark~1.2}, this was claimed also for the discrete space model with removal and without delay, though the argument is a bit trickier than it appears. 
To be precise, in the discrete space setting of \cite{Dauvergne22}*{Remark~1.2}, \eqref{eq:linear:growth} is imported from \cite{Kesten05} (see Appendix~\ref{app} for the continuous space setting), 
\eqref{eq:N:bound}, \eqref{eq:infection:number} and \eqref{eq:few:xi} are somewhat suggested to hold, while the somewhat subtle argument in \eqref{eq:def:xi} and \eqref{eq:many:xi} is omitted altogether. 
Nevertheless, it should be noted that if both time and space are discrete, the heuristics proposed in \cite{Dauvergne22}*{Remark~1.2} are essentially exhaustive.

\subsection{The model} 
\label{subsec:model}
Fix $\lambda, \alpha > 0$. Consider a Poisson Point Process $\cP$ on $\bbR$ with intensity $\lambda$ with a particle added at 0. At the beginning, the particle at position 0 is infected and all other particles are susceptible. Each particle performs an independent standard Brownian motion
$(B_p(t))_{p\in \cP,\, t\ge 0}$. Infected particles are removed at rate $\alpha$ (that is, an exponentially distributed amount of time with mean $1/\alpha$ after becoming infected) and remain removed forever. Susceptible particles become infected instantly at the first time when they are at distance at most 1 from an infected particle. As remarked in \cite{Grimmett22}*{Section~2.4}, a construction of the process informally described above can be obtained along the lines of \cite{Kesten05} (also see Appendix~\ref{app}). 

It will be convenient for us to encode the process in terms of the empirical measure of the particles and their type as follows. For any $t\ge 0$, let $\cP_t=\sum_{p\in\cP}\delta_{(B_p(t),\,\eta_p(t))}$, where $\eta_p(t)\in\cS=\{\mS,\mI,\mR\}$ is the state (susceptible/infected/removed resp.) at time $t$ of the particle starting at position $p$ in the initial condition $\cP$. We denote by $\cI_t=\sum_{p\in\cP:\,\eta_p(t)=\mI}\delta_{B_p(t)}$ the empirical measure of the infected particles at time $t$. We further define $I(t)=\cI_t(\bbR)$ as the number of infections at time $t$, as well as the leftmost and rightmost infections $L(t)=\inf\supp(\cI_t)$ and $R(t)=\sup\supp(\cI_t)$, with the convention $\inf\varnothing=\infty$ and $\sup\varnothing=-\infty$.

\subsection{The result}
Our main result establishes that, for any non-degenerate value of the parameters, infection eventually dies out almost surely. In fact, we show that the probability that the infection survives until time $T$ is exponentially small.

\begin{theorem}\label{thm:main}
There exists $c > 0$ such that, for all sufficiently large $T > 0$,
\begin{equation}
\label{eq:main}
\bbP(\sup\{t\ge 0:I(t)>0\}\ge T)\le \mathrm{e}^{-cT}.\end{equation}
In particular, almost surely there exists $T>0$ such that $I(t)=0$ for all $t\ge T$.
\end{theorem}

\begin{remark}[Extensions]
\label{rem:extensions}
While Theorem~\ref{thm:main} is stated and proved in the cleanest setting of Section~\ref{subsec:model}, it readily extends in various ways. A non-unit range of infection and diffusion coefficient can be obtained by rescaling space and time, respectively. A drift can be added to the Brownian motions by shifting space linearly with time. Furthermore, it will be clear from the proof that one can allow infection to be transmitted at finite rate which may depend on the relative position of the susceptible and infected particles, instead of instantaneously, as long as the infection range is finite. Indeed, in such models, infection spreads more slowly than in models where the infection rate is infinite and there is no removal, so Appendix~\ref{app} still applies, while the rest of the argument does not inspect the exact infection mechanism. Our technique can likely be adapted to other diffusions and sufficiently fast decaying infection rates with unbounded range (see \cite{Grimmett22}*{Theorem~3.6}), but we prefer to avoid such technical complications. For the interested reader, we point out that the only model-dependent spots in the proof are~\eqref{eq:gaussian}, which uses a very rough Gaussian computation, and Proposition~\ref{prop:linear:expansion}. The latter relies only on the stationarity of a homogeneous Poisson point process whose points perform independent Brownian motions and its result can also be recovered in a more robust way via the approach of \cite{Kesten05}, if needed.
\end{remark}

\section{Proof of Theorem~\ref{thm:main}}
Fix $\alpha,\lambda>0$. The proof proceeds as follows. First, we use the model without removal to restrict our attention to particles in an interval growing linearly with time. By standard concentration bounds this implies that only a linear number of particles is relevant. Since each infected particle typically survives a time of order 1, most of the time there are at most a bounded number of infected particles. 
Finally, we define suitable stopping times allowing us to make independent attempts to either remove all the infected particles, before they infect anyone else, or infect yet more particles. In the former case, we are done (infection dies out), while in the latter case, we arrive at a contradiction with the bound on the number of infected particles. Next, we turn to the details. Note that in the proof of Theorem~\ref{thm:main}, we often spare integer parts for the sake of readability. 

\paragraph{Linearly localised infection}
Recall that the Brownian snails model with removal at different values of $\alpha$ does not enjoy monotonicity in $\alpha$ (note that, while infected snails typically live longer for smaller values of $\alpha$, in this setting some of them may also become infected too early and fail to transmit 
the infection further). 
However, if there is no removal at all ($\alpha=0$), it holds that the set of infections at any given time contains the one for any value of $\alpha>0$. More precisely, this inclusion holds if we couple the two processes using the same Brownian motions. We may therefore use the model without removal as an upper bound. In Appendix~\ref{app}, we discuss how to adapt the proof of \cite{Beckman18} in order to obtain $c_1,C_1>0$ depending on $\lambda$ such that for any $T$ large enough
\begin{align}
\label{eq:linear:growth}
\bbP(\cE_1)
\le e^{-c_1T}\quad \text{where}\quad \cE_1
=\left\{\sup_{t\le T}\max(-L(t),R(t))\ge C_1T\right\}
\end{align}
for the model without removal (and, therefore, also for the one with removal). We have thus localised infection to an interval of linear size. We next seek to show that only linearly many particles enter this interval.
\paragraph{Linear number of infections}
Fix $T$ large enough. We next compute the expectation of the number $N$ of particles $p\in\cP$ such that $\min_{t\le T} |B_p(t)|\le C_1T$. It is at least $1+2C_1T$ (because of the particles initially present in $[-C_1T, C_1T]$) and at most
\begin{equation}
\label{eq:gaussian}
1+4C_1T+2\int_0^T\md t\int_{2C_1T}^\infty\frac{e^{-(x-C_1T)^2/(2t)}}{\sqrt{2\pi t}}\md x\le 1+4C_1T+2T\int_{C_1\sqrt{T}}^\infty \frac{e^{-x^2/2}}{\sqrt{2\pi}}\md x\le 5C_1T,\end{equation}
as $T$ is large enough. Moreover, the random variable $N-1$ has Poisson distribution with parameter $\mathbb E[N]-1$ (which is clear from the construction in Appendix~\ref{app}). We may thus apply a standard concentration result for Poisson variables to $N-1$ (see e.g.\ \cite{Alon16}*{Theorem A.1.15}) to get
\begin{align}
\label{eq:N:bound}
\bbP(\cE_2)
\le e^{-C_1T}\quad \text{where}\quad
\cE_2
=\left\{N\ge 20C_1T\right\}.
\end{align}
Thus, so far, we know that only linearly many particles may become infected.

\paragraph{Frequently bounded number of infections}
Our next goal is to show that during a linear proportion of the time, only boundedly many particles are infected. Consider a truncation of the model in which only the first $20C_1T$ infection events are allowed to occur. More precisely, the process is as above until the time when the $20C_1T$-th particle becomes infected, thereafter infection is no longer transmitted and the remaining infections are removed at rate $\alpha$ as usual. Note that on the event $\cE_1^c\cap \cE_2^c$, the truncation will not take place until time $T$. Let $I'(t)$ be the number of infected particles in the truncated process at time $t$. Then, $\int_0^\infty I'(t)\md t$ is stochastically dominated by the sum of $20C_1T$ independent exponential random variables with mean $1/\alpha$. By the exponential Markov inequality there exist $c_2,C_2>0$ depending on $\alpha$ and $C_1$ (but not on $T$) such that
\begin{align}
\label{eq:infection:number}
    \bbP(\cE_3)
    \le e^{-c_2T}\quad \text{where}\quad
    \cE_3
    =\left\{\int_0^\infty I'(t)\md t\ge C_2T\right\}.
\end{align}
Notice that, on the event $\cE^c_1\cap \cE^c_2\cap\cE^c_3$, the set $\cT=\{t\in[0,T]:I'(t)\le 2C_2\} = \{t\in[0,T]:I(t)\le 2C_2\}$ has Lebesgue measure \begin{equation}
    \label{eq:T:bound}
    \ell(\cT)\ge T/2.
\end{equation}

Fix $\varepsilon=1/(200C_1)$. We next inductively define a sequence of stopping times at which $I(t)\le 2C_2$ by setting
\begin{align}
\label{eq:def:tau}
\tau_1=\min\{t\ge 0: I(t)\le 2C_2\}\quad\text{and, for all $i\ge 1$,}\quad \tau_{i+1}=\min\{t\ge \tau_i+\varepsilon:I(t)\le 2C_2\}
\end{align}
with $\min\varnothing=\infty$.  For each $i\ge 1$, set
\begin{equation}
\label{eq:def:xi}
\xi_i=\begin{cases}1&\text{if }\tau_i<\infty\text{ and }\bbP(\exists t\in(\tau_i,\tau_i+\varepsilon),\,\lim_{\theta\to t-}I(\theta)<I(t)\mid\cP_{\tau_i})\ge1/2,\\
0&\text{otherwise}.\end{cases}
\end{equation}
In other words, we wait until there are few infections and assess whether it is likely that more infections appear during the following time interval of length $\eps$.\footnote{Note that almost surely there is no time such that a susceptible particle is infected and an infected one is removed simultaneously.} 
The key idea is to take advantage of the stopping times when infection is likely to appear in order to contradict \eqref{eq:N:bound} rather than trying to avoid this scenario.

\paragraph{Too many infections}
We first asses the stopping times $\tau_i$ such that $\xi_i=1$, that is, it is likely that a new infection will appear soon. We argue that if many such times occur, probably too many infections will arise, conflicting \eqref{eq:N:bound}. Set $i_0=0$ and, for every $j\ge 0$, 
\[i_{j+1}=\min\left\{i>i_j:\xi_i=1\right\}\]
with $\min\varnothing=\infty$. We construct a sequence $(X_j)_{j\ge 1}$ of Bernoulli random variables as follows. For every $j\ge 1$ such that $i_j<\infty$, we define 
\[X_j=\begin{cases}1&\text{if there is }t\in(\tau_{i_j},\tau_{i_j}+\varepsilon)\text{ such that }\lim_{\theta\to t-} I(\theta)<I(t),\\
0&\text{otherwise},\end{cases}\]
and extend the sequence by defining $(X_j)_{j:\, i_j = \infty}$ as a sequence of independent $\mathrm{Bernoulli}(1/2)$ random variables, which is also independent from $(X_j)_{j:\, i_j < \infty}$.
We remark that if $i_j < \infty$, then $\xi_{i_j} = 1$ and $\tau_{i_j} < \infty$, so $X_j$ is well defined. Moreover, the Markov property and \eqref{eq:def:xi} give that $\bbP(X_j=1|(X_k)_{k<j})\ge 1/2$ for all $j\ge 1$. Hence, $(X_j)_{j\ge 1}$ is stochastically dominates a sequence of i.i.d.\ Bernoulli random variables with parameter $1/2$ and a standard concentration inequality gives
\begin{align}
\label{eq:many:xi}
\bbP(\cE_4)
\le e^{-c_3T}\quad \text{where}\quad
\cE_4
=\left\{\sum_{j=1}^{50C_1T}X_j\le 20C_1T\right\}
\end{align}
for some $c_3>0$ depending only on $C_1$.

Define $M=T/(2\varepsilon) = 100C_1T$ and $\Xi=\sum_{i=1}^{M-1} \xi_i$. We show that $\bigcap_{i=1}^4\cE^c_i$ implies that $\Xi< 50C_1T$. Indeed, on the one hand, on the event $\cE_1^c\cap \cE_2^c\cap \cE_3^c$ we have $\tau_i\le T-\eps$ for all $i\le M-1$, since, by~\eqref{eq:T:bound} and~\eqref{eq:def:tau}, 
\[T/2\le \ell(\cT)\le \varepsilon(1+|\{i\in\bbN:\tau_i\le T-\varepsilon\}|).\]
On the other hand, the event $\{\Xi\ge 50C_1T\}\cap \cE^c_4$ implies that at least $20C_1T$ infections occur before time $\tau_{i_{50C_1T}}+\varepsilon\le \tau_{M-1}+\eps\le T$, which is a contradiction with $\cE^c_2$.

\paragraph{Many chances to extinguish infection}
Now that we know that $\Xi<50C_1T$, we will recover many indices $i$ such that $\xi_i=0$, so that at each corresponding $\tau_i$, we have a good chance for the infection to die out quickly. Set $i'_0=0$ and, for all integers $j\ge 0$, define
\[i_{j+1}'=\min\left\{i>i'_j:\xi_i=0\right\}\]
with $\min\varnothing=\infty$. Note that on the event $\bigcap_{i=1}^4\cE^c_i$ (implying $\Xi < 50C_1 T$) all terms of the subsequence $(i'_j)_{j=0}^{50C_1 T}$ are at most $M-1$, so finite.

Now, we construct a sequence $(Y_j)_{j\ge 1}$ of Bernoulli random variables as follows. For every positive integer $j$ such that $\max(i'_j, \tau_{i'_j})<\infty$, let 
\[Y_j=\begin{cases}
0&\text{if } I(\tau_{i'_j}+\varepsilon)=0,\\
1&\text{otherwise},\end{cases}\]
and extend the sequence by defining $(Y_j)_{j:\, \max(i'_j, \tau_{i'_j}) = \infty}$ as a sequence of independent Bernoulli random variables with parameter $p=1-(1-e^{-\alpha\varepsilon})^{2C_2}/2$, which is also independent from $(Y_j)_{j:\, \max(i'_j, \tau_{i'_j}) < \infty}$. Again, the Markov property and~\eqref{eq:def:xi} give that $\bbP(Y_j=1|(Y_k)_{k<j})\le p$ for all $j\ge 1$. Note that here, we use that for any configuration $\cP$ (not necessarily having a particle at the origin) with at most $2C_2$ infections such that the probability of infecting another particle by time $\varepsilon$ is at most $1/2$, the probability that the infection dies out completely is at least $1-p$ (it suffices for all infections to be removed within time $\varepsilon$ without modifying anything else in the graphical construction). In particular,
\begin{align}
\label{eq:few:xi}
    \bbP(\cE_5)
    \le e^{-c_4T}\quad \text{where}\quad
    \cE_5
    =\left\{\prod_{j=1}^{50C_1T}Y_j=1\right\}
\end{align}
for some $c_4>0$ depending on $\varepsilon$, $\alpha$, $C_1$ and $C_2$ but not on $T$.

We finally claim that $\bigcap_{i=1}^5\cE_i^c\subset \{I(T)=0\}$, which will conclude the proof of Theorem~\ref{thm:main} in view of \eqref{eq:linear:growth}, \eqref{eq:N:bound}, \eqref{eq:infection:number}, \eqref{eq:many:xi} and \eqref{eq:few:xi}. Assume $\bigcap_{i=1}^5\cE_i^c$ occurs. As already discussed, $\Xi<50C_1T$, so $i'_{50C_1T}\le M-1$. But then $\cE^c_5$ yields $I(\tau_{i'_{50C_1T}}+\varepsilon)=0$, which, together with the fact that $\tau_{M-1}\le T-\eps$, implies that $I(T) = 0$, as desired. The almost sure extinction follows from~\eqref{eq:main} and the Borel--Cantelli lemma.

\section*{Acknowledgements}
This work was supported by the Austrian Science Fund (FWF): P35428-N. We thank Geoffrey Grimmett for introducing us to the problem and for encouraging remarks. We are also grateful to the organisers of the Recent Developments in Stochastic Processes conference, which sparked this project.

\appendix
\section{At most linear growth}
\label{app}
In this appendix, we discuss the proof of \eqref{eq:linear:growth}. Since the argument is quite general, we work directly in arbitrary dimension $d\ge 1$.

Let us begin by introducing the Brownian snails model without removal (that is, with $\alpha=0$) more formally. Fix $\lambda\in(0,\infty)$. Let $\cP$ be a Poisson point process on $\cC_0\times \bbR^{d}$, where $\cC_0$ denotes the space of continuous functions $f:[0,\infty)\to\bbR^d$ such that $f(0)=0$ equipped with the topology of uniform convergence on compact sets (see e.g.\ \cite{Kingman93}*{Chap.~2} for background on Poisson processes). We refer to elements of $\cC_0\times\bbR^d$ as \emph{possible trajectories} and to elements of $\cP$ (or other Poisson point processes on the same space) as \emph{trajectories}. We view the second coordinate of a possible trajectory as its \emph{starting point}. We take the intensity measure of $\cP$ to be $\bbW_d\otimes(\lambda\ell)$, where $\bbW_d$ denotes the $d$-dimensional Wiener measure (which is the distribution of the standard Brownian motion on $\bbR^d$ started at 0) and $\ell$ is the Lebesgue measure on $\bbR^{d}$. We further fix an independent $d$-dimensional Brownian motion $B_0$ and set $\cP_0=(\cP\setminus(\cC_0\times \bbB(0,1)))\cup\{(B_0,0)\}$, where $\bbB(0,1)$ is the unit ball of $\bbR^d$. Thus, one obtains $\cP_0$ from $\cP$ by removing the trajectories starting within unit distance from the origin, and then adding one trajectory starting at the origin. We label the trajectories in $\cP_0$ by $((B_i(t))_{t\ge 0},x_i)_{i=0}^\infty$. It remains to specify the state ($\mS$ or $\mI$) of each particle. We set $T(0)=0$ and for any positive integer $i$, we define 
\begin{multline*}
T(i)=\inf\{t>0:\exists k\ge 0,\,\exists i_0=0,\exists i_1,\dots,\exists i_k=i,\,\exists 0<t_1\le t_2\le \dots \le t_k\le t,\,\forall j\in\{1,\dots,k\},\\\|(B_{i_j}(t_j)+x_{i_j})-(B_{i_{j-1}}(t_j)+x_{i_{j-1}})\|\le 1\}.\end{multline*}
Thus, the positions 
of the infected particles at time $t\ge 0$ are 
given by 
\[I_t=\left\{B_i(t)+x_i:i\ge 0,\, T(i)\le t\right\}.\]
With this notation, we are ready to state the result we are after.
\begin{proposition}
\label{prop:linear:expansion}
    Let $d\ge 1$. For any $\lambda>0$ strictly smaller than the critical rate for continuum percolation with radius $1$ in $d$ dimensions (equal to $\infty$ when $d=1$; see \cite{Penrose03}), there exist $c,C>0$ such that, for all $t$ large enough,
    \[\bbP\left(\max\left\{\|z\|:z\in I_t\right\}\ge Ct\right)\le e^{-ct}.\]
\end{proposition}
While one can probably use the method of \cite{Kesten05}*{Theorem~1} (with a significant amount of work, in particular due to the need of also discretising space and not only time) to prove this fact, and it was more or less claimed in \cite{Grimmett22}*{Section 2.4}, there is a more elegant approach. It turns out that, once one adopts the right viewpoint, the proof becomes essentially identical to the one of \cite{Beckman18}. Indeed, \cite{Beckman18}*{Proposition 1.4} is the analogue of Proposition~\ref{prop:linear:expansion} for the corresponding model with delay (where particles remain immobile until they are infected). The gist of \cite{Beckman18} in our setting is as follows. 

We start by exploring $B_0$ until the time $\tau=\inf_{i\ge 1}T(i)$ when the first infection occurs. Also, reveal the (random) index $\iota$ of the particle which comes at distance 1 from $B_0(\tau)$ at time $\tau$, together with its position $B_\iota(\tau)+x_\iota$ at time $\tau$. Then comes the key point: one considers the random set \[X=\left\{(B,x)\in \cC_0\times\bbR^d:\exists t\le \tau,\, \left\|B(t)+x-B_0(t)\right\|\le 1\right\}\]
of possible trajectories coming within distance $1$ of $B_0$ no later than time $\tau$. We know by construction that $\cP_0\cap X=\{(B_0,0),(B_{\iota},x_{\iota})\}$. Therefore, $X$ is measurable with respect to (the sigma-algebra generated by) $\cP_0\cap X$, since 
\[\tau=\min\{t\ge 0:\|B_\iota(t)+x_\iota-B_{0}(t)\|=1\}\]
is measurable with respect to $\cP_0\cap X$. 

We then define a Poisson point process $\cP'$ that is equal to $\cP_0$ outside $X$ and independent with intensity $\1_X(\bbW_d\otimes (\lambda\ell))$ in $X$. Observe that $\cP'$ is indeed a Poisson point process equal to $\cP$ in distribution. To see this, it suffices to note that for any disjoint measurable $Y_1,\dots, Y_k\subset \cC_0\times \bbR^d$, the random variables
\begin{equation}
\label{eq:decomposition}
|\cP'\cap Y_i|=|(\cP_0\setminus X)\cap Y_i|+|\cP'\cap Y_i\cap X|\end{equation}
are independent for all $i\in\{1,\dots,k\}$ and have expectations $(\bbW_d\otimes(\lambda\ell))(Y_i)$. Moreover, in view of \eqref{eq:decomposition}, $\cP'$ is independent of $X$, since $X$ is measurable with respect to $\cP_0\cap X$. Let $\cP'=\{(B'_{i},x'_i):i\ge 0\}$ and write $(B'_{-1},x'_{-1})=(B_\iota,x_\iota)$ and $(B'_{-2},x'_{-2})=(B_0,0)$ for convenience.

We then consider the continuum percolation cluster 
\begin{multline*}
C=\Big\{j\ge 0:\exists k\ge0,\exists i_0\ge 0,\exists i_1\ge 0,\dots,\exists i_k=j,\;\left\|B'_{i_0}(\tau)+x'_{i_0}-B_\iota(\tau)-x_\iota\right\|\le 1\text{ and }\\
\forall j'\in\{1,\dots,k\},\,
\left\|B'_{i_{j'}}(\tau)+x'_{i_{j'}}-B'_{i_{j'-1}}(\tau)-x'_{i_{j'-1}}\right\|\le 1\Big\}
\end{multline*}
of the point $B_{\iota}(\tau)+x_{\iota}$ with radius 1 in the projection on $\bbR^d$ at time $\tau$ of the trajectories in $\cP'$ (that is, the set of particles other than $\iota$ that would become immediately infected if we introduce an infection at $B_{\iota}(\tau)+x_{\iota}$ at time $\tau$). Finally, for each $j\in C\cup\{-2,-1\}$, we start a copy of the original process translated in time by $\tau$ and in space to $B'_j(\tau)+x'_j$ (note that here we are using the fact that for every $t\ge 0$, $\{B(t)+x:(B,x)\in\cP\}$ is still a Poisson point process on $\bbR^d$ with intensity $\lambda \ell$). The Poisson point process for each $j\in C\cup\{-2,-1\}$ is given by $(\cP'\setminus\{(B'_{j'
},x'_{j'
}):j'\in C\})\cup\{(B'_{j},x'_j)\}$ to which we add an independent Poisson point process with intensity $\bbW_d\otimes (\lambda\ell)$ restricted to the set 
\[\left\{(B,x)\in\cC_0\times\bbR^d: \left\|B(\tau)+x-B'_j(\tau)-x'_j\right\|>1,\exists j'\in C\cup\{-2,-1\},\,\left\|B(\tau)+x-B'_{j'}(\tau)-x'_{j'}\right\|\le 1\right\}\]
of possible trajectories which at time $\tau$ are located at distance less than 1 from $B'_{j'}(\tau)+x'_{j'}$ for some $j'\in C\cup\{-2,-1\}$, but at distance more than 1 from $B'_{j}(\tau)+x_{j}$.

One then studies the branching process obtained this way and proves that it grows at most linearly. Since the only difference with respect to \cite{Beckman18}*{Section~3 and Appendix~A} is the use of the space $\cC_0\times \bbR^{d}$ instead of $\bbR^d$, we direct the reader to that paper for more details. 
The only subtlety to account for is the fact that $\bbP(\tau=0)=0$, which is slightly more difficult when all particles move, but is proved as in \eqref{eq:N:bound}.

\begin{remark}
A similar approach can be used to obtain a somewhat simpler proof of \cite{Kesten05}*{Theorem 1}.
\end{remark}
Finally, let us note that, while Proposition~\ref{prop:linear:expansion} only gives a bound at a given time, it is not hard to deduce that $\bbP(\max\{\|z\|:\,\exists t'\le t,\, z\in I_{t'}\}\ge 2Ct)\le e^{-c't}$ for some $c'>0$, as claimed in \eqref{eq:linear:growth}. Indeed, by a reasoning similar to~\eqref{eq:N:bound}, if an infected particle is at distance $2Ct$ from the origin at some time $t'\le t$, it is exponentially unlikely to be within distance $Ct$ of the origin at time $t$.

\bibliographystyle{plain}
\bibliography{Bib}

\end{document}